\documentclass[a4paper,reqno]{amsart}
\usepackage{amssymb}
\usepackage{latexsym}
\usepackage{amsmath}
\usepackage{euscript}
\usepackage{graphics,color}
\usepackage[all]{xy}
\usepackage[margin=3cm]{geometry}
\usepackage{MnSymbol} 

\newcommand{\be}{\begin{equation}}
\newcommand{\ee}{\end{equation}}
\newcommand{\ba}{\begin{eqnarray}}
\newcommand{\ea}{\end{eqnarray}}
\newcommand{\baa}{\begin{eqnarray*}}
\newcommand{\eaa}{\end{eqnarray*}}
\newcommand{\bb}{\begin{thebibliography}}
\newcommand{\eb}{\end{thebibliography}}

\newcommand{\bi}[1]{\bibitem{#1}}
\newcommand{\lab}[1]{\label{#1}}
\newcommand{\re}[1]{(\ref{#1})}

\newcounter{my}
\newcommand{\he}%
   {\stepcounter{equation}\setcounter{my}%
   {\value{equation}}\setcounter{equation}0%
   }%
\newcommand{\she}%
   {\setcounter{equation}{\value{my}}%
    }%

\renewcommand\t{\tilde}

\newcommand\vphi{\varphi}

\newtheorem{pr}{Proposition}

\theoremstyle{definition}

\numberwithin{equation}{section}

\begin{document}

\title[OPUC Kronecker polynomials]{Para-orthogonal polynomials on the unit circle generated by Kronecker polynomials}



\author{Alexei Zhedanov}

\address{School of Mathematics, Renmin University of China, Beijing 100872, China}

\vspace*{5mm}

\begin{abstract}
The Kronecker polynomial $K(z)$ is a finite product of cyclotomic polynomials $C_j(z)$. Any Kronecker polynomial $K(z)$ of degree $N+1$ with simple roots  on the unit circle generates a finite set $\Phi_0=1, \Phi_1(z), \dots, \Phi_N(z) $ of polynomials (para) orthogonal on the unit circle (POPUC).  This set is determined uniquely by the condition $\Phi_N(z) = (N+1)^{-1} K'(z)$. Such set can be called the set of Sturmian Kronecker POPUC.  We present several new explicit examples of such POPUC. In particular, we define and analyze properties of the Sturmian cyclotomic POPUC generated by the cyclotomic polynomials $C_M(z)$. Expressions of these polynomials strongly depend on the decomposition of $M$ into prime factors.   
 \end{abstract}

\maketitle

\section{Introduction}
\setcounter{equation}{0} 
Assume that $P(x)$ is an arbitrary monic polynomial (i.e. with leading term equal to 1) with all real simple zeros:
\be
P_{N+1}(x) = (x-x_0)(x-x_1) \dots (x-x_N), \quad x_0<x_1<x_2<\dots < x_N \lab{P_pr} \ee 
The derivative polynomial
\be
P_N(x) = (N+1)^{-1} P_{N+1}'(x)
\lab{P_N} \ee
has all zeros real and simple as well. Moreover, the zeros of $P_N(x)$ interlace zeros of $P_{N+1}(x)$. This means that each interval $(x_i,x_{i+1})$ between two neighbor zeros of $P_{N+1}(x)$ contains one and only one zero of $P_N(x)$.

Then the classical Sturm algorithm \cite{Waerden} generates the finite sequence of orthogonal polynomials $P_{N-1}(x), P_{N-2}(x), \dots , P_0=1$ such that $\deg (P_n(x)) =n$. The three-term recurrence relation between these polynomials 
\be
P_{n+1}(x) +(b_n-x) P_n(x) + u_n P_{n-1}(x)=0
\lab{3term_P} \ee 
is a consequence of the Euclidean division algorithm with respect to the pair of polynomials $P_{n+1}(x), P_n(x)$ (the next polynomial $P_{n-1}(x)$ appears as the remainder).  Note that necessarily $u_n>0$ for all $n=1,2,\dots, N$. This condition follows easily from the initial condition that all zeros of $P_{N+1}(x)$ are simple and real.  

Moreover, it is elementary  to show that the polynomials $P_n(x)$ satisfy the orthogonality relation \cite{Chi}
\be
\sum_{s=0}^N P_n(x_s) P_m(x_s) w_s = h_n, \delta_{nm} \lab{ort_P} \ee 
where the weights are given by the expression \cite{ZheS}
\be
w_s = \frac{h_N}{{P_{N+1}'(x_s)}^2} \lab{w_s_P} \ee
and where 
\be
h_N=u_1 u_2 \dots u_N >0 \lab{h_N} \ee
is the normalization factor. Note that the weights are normalized by the standard condition
\be
\sum_{s=1}^N w_s =1. \lab{norm_w} \ee

It is natural to refer the sequence $P_{N+1}(x), P_N(x), \dots P_0(x)$ as the {\it Sturmian sequence} of orthogonal polynomials \cite{ZheS}. 

\vspace{3mm}

In \cite{ZheS} it was demonstrated how some "classical" Sturmian sequences can be generated from a polynomial $P_{N+1}(x)$ with prescribed zeros on linear or quadratic grids.  The purpose of this paper is to consider  analogs and applications of the Sturmian sequence for the case of polynomials orthogonal on the unit circle (OPUC for brevity). Namely, we apply corresponding algorithm to the special case of the {\it Kronecker polynomials} in order to obtain many new explicit examples of OPUC.

We first recall basic facts on OPUC (see \cite{Simon} for details).  

Let $\Phi_n(z)=z^n + O(z^{n-1})$ be a set of monic polynomials satisfying the recurrence relation \cite{Simon}
\be
\Phi_{n+1}(z) =z \Phi_n(z) - \bar a_n \Phi^*_n(z), \lab{rec_Phi} \ee 
where 
$$
\Phi_n^*(z) =z^n \bar \Phi_n(1/z)
$$
and where $\bar \Phi_n(z)$ means the polynomials obtained by complex conjugation of expansion coefficients of the polynomial $\Phi_n(z)$ . The parameters
$$
a_n = -\bar \Phi_{n+1}(0)
$$ 
play the crucial role in the theory and properties of OPUC (these parameters are called sometimes the Verblunsky parameters \cite{Simon}).

Relation \re{rec_Phi} can be presented in an equivalent form 
\be
\Phi^*_{n+1}(z) = \Phi^*_{n}(z) - z a_n \Phi_n(z) \lab{rec_Phi*} \ee

One of the main result in the theory of OPUC states that if  
\be
|a_n|<1, \quad n=0,1,2,\dots , \lab{a_cond} \ee
then the polynomials $\Phi_n(z)$ are orthogonal on the unit circle with respect to a positive measure $d \sigma(\theta)$
\be
\int_0^{2 \pi} \Phi_n(e^{i \theta}) \bar \Phi_m(e^{-i\theta}) d \sigma(\theta) = h_n \delta_{nm} \lab{ort_Phi} \ee
where
\be
h_n = (1-|a_0|^2)(1-|a_1|^2) \dots (1-|a_{n-1}|^2) \lab{h_n} \ee
are nonzero (in fact, positive) normalization constants.

Note that orthogonality relation\re{ort_Phi} is equivalent to conditions \cite{Simon}
\ba
\int_0^{2 \pi} \Phi_n(e^{i \theta}) e^{-ij\theta} d \sigma(\theta) = h_n \delta_{nj}, \quad j=0,1,2,\dots, n  \lab{ort_j} \ea

Moreover, it can be showed that all zeros of polynomials $\Phi_n(z)$ lie inside the unit circle $|z|<1$. In contrast to the case of polynomials $P_n(x)$ orthogonal on the real axis, the roots of $\Phi_n(z)$ need not be simple \cite{Simon}.

There is an important {\it inverse} formula which allows to find the polynomial $\Phi_n(z)$ if the polynomial $\Phi_{n+1}(z)$ is known. This relation follows from relations \re{rec_Phi} - \re{rec_Phi*} \cite{Simon}
\be
\Phi_n(z) = \frac{\Phi_{n+1}(z) + \bar a_n \Phi^*_{n+1}(z)}{z(1-|a_n|^2)}
\lab{inv_Phi} \ee
Relation \re{inv_Phi} uniquely determines $\Phi_n(z)$ from the given polynomial $\Phi_{n+1}(z)$ if $|a_n| <1$.

The case $|a_n|=1$ is exceptional but it is important because it leads to  a finite system of OPUC sometimes called the {\it para-orthogonal} polynomials (POPUC) \cite{MS}, \cite{Mar1} .

Indeed, assume that $|a_i|<1, \: i=0,1,\dots, N-1$ but that $|a_N|=1$. Then from recurrence relation \re{rec_Phi} it follows that the polynomial $\Phi_{N+1}(z)$ has distinct roots at the unit circle:
\be
\Phi_{N+1}(z) = (z- \zeta_0)(z-\zeta_1) \dots (z-\zeta_N) \lab{Phi_unit} \ee
where
\be
|\zeta_i| =1, \; i=0,1,\dots, N.\lab{zeta_unit} \ee  
It is easy to show that in this case the polynomials $\Phi_k(z), \, k=0,1,\dots, N$ are orthogonal on the unit circle with discrete measure located at the zeros of the polynomial $\Phi_{N+1}(z)$
\be
\sum_{s=0}^N w_s \Phi_n(\zeta_s) {\bar \Phi}_m(\zeta^{-1}_s) = h_n \delta_{nm},
\lab{ort_fin} \ee
where the positive weights are \cite{Ger}, \cite{Simon}
\be
w_s = \frac{h_N}{\Phi_{N+1}'(\zeta_s) {\bar \Phi}_N(\zeta^{-1}_s)}.
\lab{w_s_OPUC} \ee

Explicit examples of polynomials orthogonal on unit circle are very interesting from different points view.  A list of known explicit examples can be found e.g. in Simon's monograph \cite{Simon}.

In this paper we propose new explicit examples of finite systems POPUC with concentrated masses located at roots of unity. More exactly, we define the system $\Phi_0(z)=1, \Phi_1(z), \dots, \Phi_{N}(z), \Phi_{N+1}(z)$ by the conditions
\be
\Phi_{N+1}(z) = K_{N+1}(z), \quad \Phi_{N}(z) = (N+1)^{-1} K'_{N+1}(z),
\lab{Sturm_pair} \ee 
where $K_{N+1}(z)$ is the Kronecker polynomial with $N+1$ simple roots on the unit circle. 

The paper is organized as follows. In Section 2 we describe the Sturmian algorithm to construct a finite set of POPUC starting from prescribed polynomial $\Phi_{N+1}(z)$. In Section 3, we recall definition and properties of the Kronecker polynomials $K(z)$ and define the Sturmian Kronecker POPUC. In Section 4 we consider the simplest choices of $K(z)$ leading to elementary expressions of Sturmian Kronecker POPUC. In Section 5,  the properties of the Sturmian cyclotomic POPUC  are considered. These POPUC are generated by a single cyclotomic polynomial $C_M(z)$. In Section 6 some explicit examples of Sturmian Kronecker POPUC beyond the cyclotomic polynomials are presented. Finally, in Conclusion, we discuss some open problem.

\section{Sturmian sequence for the unit circle}
\setcounter{equation}{0}
In this section we follow basically the approach proposed in \cite{CMR} with some modifications.  Assume that a polynomial $\Phi_{N+1}(z)$ is given with all simple zeros on the unit circle:
\be
\Phi_{N+1}(z) = (z-\zeta_0) (z-\zeta_1) \dots (z-\zeta_N),
\lab{Phi_zeta} \ee
where $|\zeta_i|=1, \: i=0,1,\dots, N$.
We assume also that there exists a sequence of OPUC $\Phi_i(z), \: i=0,1,\dots, N$ such that $\Phi_{N+1}(z)$ is the final member of this sequence. This means in particular that the relation
\be
\Phi_{N+1}(z) = z \Phi_N(z) - \bar a_N \Phi^*_{N}(z) \lab{rel_NN} \ee 
should be valid, where
\be
\bar a_N = - \Phi_{N+1}(0) = (-1)^N \zeta_0 \zeta_1 \dots \zeta_N. \lab{a_N_zeta} \ee 
Then starting with $\Phi_N(z)$ we can construct uniquely the next monic polynomial  $\Phi_{N-1}(z)$ by using the inverse formula \re{inv_Phi}. Obviously this process can be continued to produce $\Phi_{N-2}(z), \Phi_{N-3}(z), \dots$ until achieving the last member $\Phi_0=1$.

Thus the pair of polynomials $\Phi_{N+1}(z)$ and $\Phi_N(z)$ completely determines the whole finite sequence $\Phi_i(z), \: i=0,1,\dots, N$ of OPUC which is orthogonal on zeros $\zeta_s$ of the polynomial $\Phi_{N+1}(z)$.

If only $\Phi_{N+1}(z)$ is given, how general could be the polynomial 
$\Phi_N(z)$? Obviously, it is sufficient to check the restrictions on 
$\Phi_N(z)$ followed from relation \re{rel_NN}. This leads to {\it infinitely many} possible candidates for the appropriate polynomial $\Phi_N(z)$.

All these candidates can be explicitly presented using the Wendroff theorem \cite{Mar1} which states that the set of OPUC $\Phi_0,\Phi_1(z), \dots, \Phi_N(z)$ is uniquely determined by two given sets $z_0, z_1, \dots z_{N}$ and $\t z_0, \t z_1, \dots \t z_{N}$ of points on the unit circle (i.e. $|z_i|=|\t z_i|=1$ for all $i=0,1,\dots, N$) with the only interlacing property which means that each arc connecting two neighbor points from the set $\{ z_i \}$ contains one and only one point from the set $\{\t z_i\}$. 

Let $\Phi_{N+1}(z) = (z-z_0)(z-z_1) \dots (z-z_N)$ and $\t \Phi_{N+1}(z) = (z-\t z_0)(z- \t z_1) \dots (z-\t z_N)$ be corresponding characteristic polynomials of these two sets. Then one has the two relations
\be
\Phi_{N+1}(z) = z \Phi_N(z) - \bar a_N \Phi_N^*(z) \quad \mbox{and} \quad \t \Phi_{N+1}(z) = z  \Phi_N(z) - {\t {\bar a}}_N \Phi_N^*(z), \lab{PPP} \ee
where
\be
a_N = (-1)^{N+1} {\bar z_0} {\bar z_1} \dots {\bar z_N}, \quad \t a_N = (-1)^{N+1}  {\t {\bar z}}_0  {\t {\bar z}}_1 \dots  {\t {\bar z}}_N , \quad |a_N| =|\t a_N| =1. \lab{aa} \ee
From \re{PPP} we obtain the expression of the polynomial
\be
\Phi_{N}(z) = \frac{{\t {\bar a}}_N \Phi_{N+1}(z) - \bar a_N \t \Phi_{N+1}(z)}{({\t {\bar a}}_N - \bar a_N)z}.
\lab{Phi_N_W} \ee
By the Wendroff theorem, the polynomial $\Phi_N(z)$ in \re{Phi_N_W} has all roots inside the unit circle $|z|<1$ (these roots are not necessarily distinct).

There is one special choice of $\Phi_N(z)$ which corresponds to the ordinary Sturm algorithm.

Indeed, let us consider the monic polynomial $\Phi_N(z)$ of degree $N$ as the {\it derivative} of the initial polynomial $\Phi_{N+1}(z)$:  $\Phi_{N}(z) = (N+1)^{-1} \Phi'_{N+1}(z)$. By the Gauss-Lucas theorem, all roots of the polynomial $\Phi_{N}(z)$ belong to the open unit disc $|z|<1$.

Moreover, it is easily seen that relation \re{rel_NN} holds. Indeed, the polynomial $\Phi_{N+1}(z)$ satisfy the obvious relation
\be
\Phi^*_{N+1}(z) = -a_N \Phi_{N+1}(z). \lab{N+1_sym} \ee 
Taking derivative of \re{N+1_sym} we arrive at relation \re{rel_NN}. This means that the polynomial $\Phi_N(z) = (N+1)^{-1}\Phi'_{N+1}(z)$ is admissible as the next member of the POPUC sequence.

We thus see that similarly to the real case, the Sturmian pair $\Phi_{N+1}(z)$ and $(N+1)^{-1}\Phi'_{N+1}(z)$ generates a unique finite sequence of POPUC. This is one of the main results obtained  in \cite{CMR}.

These polynomials are orthogonal 
\be
\sum_{s=0}^N w_s \Phi_n(\zeta_s) {\bar \Phi}_m(\zeta_s^{-1}) = h_n \delta_{nm} 
\lab{ort_Phi_St} \ee
where
\be
w_s = \frac{h_N}{\left|\Phi'_{N+1}(\zeta_s)\right|^2}
\lab{ws_st} \ee
are positive weights. Formula \re{ws_st} follows directly from \re{w_s_OPUC}.

Hence the problem is reduced to an appropriate choice of the polynomial $\Phi_{N+1}(z)$ having all its $N+1$ simple roots on the unit circle. 

We propose here a method of the {\it inverse Sturm problem}. Recall that for the case of the polynomials orthogonal on the real line, the Sturm algorithm allows to determine location of zeros of the polynomial $P_{N+1}(x)$ if these zeros are unknown initially \cite{Waerden}. In \cite{ZheS} the inverse Sturm problem was proposed. This means that we start with a polynomial $P_{N+1}(x)$ with {\it prescribed} simple zeros on the real line. Then the Strum algorithm allows to reconstruct the whole chain $P_{N-1}(x), P_{N-2}(x), \dots, P_0(x)$ of polynomials orthogonal on these prescribed zeros. It appears that for zeros of $P_{N+1}(x)$ belonging to some "classical" grids (e.g. linear or quadratic) one can generate some classical (or semiclassical) systems of orthogonal polynomials (see \cite{ZheS} for details).

Similarly, we can start with a polynomial $\Phi_{N+1}(z)$ with {\it prescribed zeros} on the unit circle.  The polynomial $\Phi_N(z)$ is chosen as the derivative of $\Phi_{N+1}(z)$. Then the polynomials $\Phi_{N-1}(z), \Phi_{N-2}(z), \dots, \Phi_0$ are reconstructed uniquely via inverse formula \re{inv_Phi}. It is naturally to call the sequence $\Phi_{N+1}(z), \Phi_N(z), \Phi_{N-1}(z),  \dots, \Phi_0=1$ the {\it Strmian sequence of POPUC}.

The main problem here is an appropriate choice of the polynomial $P_{N+1}(z) $ with prescribed zeros. The simplest and natural choice is to take zeros of the polynomial $\Phi_{N+1}(z)$ coinciding with sets (or subsets) of roots of unity. Roughly speaking, the set of roots of unity can be considered as a circle analog of the linear grid on the real axis. The most natural choice of the polynomial $\Phi_{N+1}(z)$ for this purpose is to take the {\it Kronecker polynomials}. In the next section we recall basic properties of these polynomials.

It should be stressed that the Kronecker polynomials are not the only polynomials with prescribed zeros at roots of unity. There are infinitely many such polynomials beyond the set of Kronecker polynomials.  We have decided to restrict ourselves with the Kronecker polynomials only because of their nice properties leading to new explicit examples of POPUC.

\section{Kronecker polynomials and their Sturmian POPUC}
\setcounter{equation}{0}
It is naturally to choose the polynomial $\Phi_{N+1}(z)$ coinciding with the {\it Kronecker polynomial} (this name was proposed in \cite{Dam}). We recall here the definition and basic properties of these polynomials. See \cite{Dam} for further details.  

In 1857 Leopold Kronecker introduced \cite{Kron} a set of monic polynomials with integer coefficients
\be
K(z) = z^{N+1} + a_N z^N + \dots a_k z^k + \dots + a_0, \quad a_k \in \mathbb{Z}.
\lab{KP} \ee  
The problem which Kronecker had proposed and solved was to describe effectively all such polynomials with the additional restriction: all their roots are assumed to lie in  the (closed) unit disk $|z|\le 1$.

In what follows we assume that $a_0 \ne 0$. Then the polynomial $K(z)$ has all nonzero roots. If nevertheless $a_0=0$, then $K(z)$ has a (probably multiple) zero root. In this case we can introduce the new Kronecker polynomial $\t K(z)$ by the relation
\be
K(z) = z^j \t K(z), 
\lab{KtK} \ee 
where $j$ is the multiplicity of the zero root. Obviously the polynomial $\t K(z)$ has no zero roots. This observation allows to restrict consideration only with the polynomials $K(z)$ with $a_0  \ne 0$ (i.e. with nonzero roots).

Kronecker established the following properties of the polynomials $K(z)$:

(i) all roots $\zeta_i$ of the polynomial $K(z)$ lie on the unit circle: $|\zeta_i| =1$. Moreover, every root $\zeta_i$ is a root of unity.

\vspace{3mm}

(ii) for the given $N=1,2,\dots$ there is only a finite number of distinct Kronecker polynomials of degree $N$.

\vspace{3mm}

(iii) every Kronecker polynomial $K(z)$ can be presented as a finite product of {\it cyclotomic polynomials}
\be
K(z) = C_{m_1}^{j_1}(z) C_{m_2}^{j_2}(z) \dots C_{m_k}^{j_k}(z) 
\lab{K_C} \ee  
with positive integers $m_i$ and $j_i$. 

Recall that the cyclotomic polynomial $C_n(z)$ is the minimal polynomial of all primitive $n$-th roots of unity. The primitive $n$-th root  of unity $\zeta$ is defined as a complex number such that $\zeta^n=1$ but $\zeta^k \ne 1$ for $k=1,2, \dots, n-1$. It is well known that for given $n$ the total number of primitive roots of unity is $\vphi(n)$, where $\vphi(n)$ is the Euler totient function (i.e. the number of all positive integers smaller than $n$ and coprime with $n$).

If $n=p$ is a prime number,  then all roots of unity (apart from the trivial $z=1$) are primitive and hence the cyclotomic polynomial
\be
C_p(z) = \frac{z^p-1}{z-1} = z^{p-1} + z^{p-2} + \dots + z+1
\lab{C_prime} \ee 
has degree $p-1$. For non-prime integers $n$ the explicit expression of the cyclotomic polynomial $C_n(z)$ strongly depends on decomposition of $n$ into prime factors.

For example, for $n=6$ there are only two primitive roots of unity and hence degree of the cyclotomic polynomial is two:
\be
C_6(z) = z^2-z+1 \lab{C_6} \ee 
For $n=8$ there are  four primitive roots of unity and
\be
C_8(z) = z^4+1.
\lab{C_8} \ee
The cyclotomic polynomials possess many remarkable properties, e.g. they are monic polynomials with integer coefficients irreducible over the field of the rational numbers. Moreover there is the relation
\be
x^n-1 = \prod_{d|n} C_d(z),
\lab{prod_C} \ee
where $d$ runs over all  positive divisors of $n$. Formula \re{prod_C} allows to express cyclotomic polynomials inductively using expressions of smaller polynomials.

Thus for any Kronecker polynomial (without zero roots) all its roots belong to a set of roots of unity. This property makes these polynomials very convenient for the choice of the "final" polynomial $\Phi_{N+1}(z)$ of a finite POPUC set $\Phi_0, \Phi_1(z), \dots, \Phi_N(z), \Phi_{N+1}(z)$. 

Namely, our main goal will be construction of the following Sturmian sequence of POPUC. We start with any Kronecker polynomial $K(z)$ of degree $N+1$ with simple nonzero roots and then define  
\be \Phi_{N+1}(z) = K(z), \quad \Phi_N(z) =\frac{K'(z)}{N+1}
\lab{start_K} \ee
By the above considerations, the polynomial $\Phi_{N}(z)$ has all roots inside the unit circle. Hence we can define the whole Sturmian sequence of monic polynomials $\Phi_{N-1}(z), \Phi_{N-2}, \dots, \Phi_0(z)=1$. These polynomials are unique and satisfy recurrence relation  \re{rec_Phi} with $|a_n|<1$. Hence these polynomials constitute a set of POPUC. We will call these polynomials the {\it Sturmian Kronecker POPUC}.

We notice the first simple (almost obvious) properties of these polynomials

\begin{pr}
Any admissible Kronecker polynomial $K(z)$ (i.e. which generates the set of Sturmian Kronecker POPUC $\Phi_n(z), \: n=0,1,\dots, N$) has the presentation 
\be
K(z) = C_{m_1}(z) C_{m_2}(z) \dots C_{m_k}(z), 
\lab{K_adm} \ee
where $C_{m_i}(z)$ are distinct cyclotomic polynomials.
\end{pr}
Indeed, formula \re{K_C} describes an arbitrary Kronecker polynomial without zero roots. In order to obtain the Sturmian Kronecker set of POPUC we should demand that all roots of the polynomial $K(z)$ be distinct. This leads to the condition $j_1=j_2=\dots = j_k=1$ in \re{K_C}. It remains to note that the distinct cyclotomic polynomials have distinct roots. This leads to formula \re{K_adm}.

\begin{pr}
The expansion coefficients $A_{ni}$ of Kronecker POPUC $\Phi_n(z)=\sum_{i=0}^n A_{ni} z^i$ are rational numbers.
\end{pr}
The proof of this proposition is almost obvious. We first observe that the polynomial $\Phi_N(z)$ defined by \re{start_K},  has rational expansion coefficients. Then all further polynomials  $\Phi_{N-1}(z), \Phi_{N-2}, \dots, \Phi_0(z)=1$ have rational expansion coefficients as well. This follows immediately from formula \re{inv_Phi}.

In particular, all Verblunsky coefficients $a_n$ are rational numbers.

In the next section we consider two simplest cases of the Kronecker polynomials and corresponding systems of POPUC.

\section{Simplest examples of Sturmian Kronecker POPUC}
\setcounter{equation}{0}
In order to construct a Sturmian sequence of POPUC we can start with the Kronecker polynomial $K(z)$ of degree $N+1$ with all  $N+1$-th roots of unity on the unit circle. Then as we know, the polynomial $\Phi_N(z) = (N+1)^{-1} K'(z)$ has all its roots inside of the unit circle. The inverse algorithm \re{inv_Phi} will then uniquely determine the whole chain $\Phi_{N-1}(z), \Phi_{N-2}(z), \dots \Phi_0(z)=1$ of POPUC.

Obviously, this polynomial is 
\be
K_{N+1}(z) = z^{N+1} -1.
\lab{sK} \ee
Its zeros are
\be
\zeta_k = \exp\left( \frac{2 \pi i k}{N+1} \right), \quad k=0,1,\dots, N
\lab{zeta_s} \ee
Clearly,
\be
\Phi_{N}(z) = (N+1)^{-1} K_{N+1}'(z) = z^N.
\lab{Phi_Ns} \ee
By induction we have 
\be
\Phi_k(z) = z^k, \quad k=0,1,2,\dots, N.
\lab{Phik_s} \ee
All Verblusnky parameters (apart from $a_N$) are zero 
\be
a_0=a_1=\dots =a_{N-1}=0, \; a_N=1.
\lab{a_s} \ee
This corresponds to so-called "free" POPUC \cite{Simon}. They are orthogonal on vertices $z=\zeta_s$ of the regular $N+1$-gon with equal concentrated masses
\be
\sum_{s=0}^{N} \Phi_n(\zeta_s) \Phi_m(\zeta_s^{-1}) =(N+1) \delta_{nm}.
\lab{ort_s} \ee

Consider a less trivial case with the Kronecker polynomial  
\be
K(z)= \Phi_{N+1}= \frac{z^{N+2}-1}{z-1} = z^{N+1} + z^N + \dots + z+1 \lab{1free_N+1} \ee
The roots are
\be
\zeta_k = \exp\left(\frac{2\pi i k}{N+2} \right), \quad k=1,2,\dots, N+1 \lab{1zeta_reg} \ee
I.e. the roots are located at vertices of regular $N+2$-gon apart from the point $z=1$.
We have
\be
\Phi_N(z) = \frac{K'(z)}{N+1} = \frac{(N+1) z^{N} + N z^{N-1} + \dots + k z^{k-1} + 2 z + 1}{N+1} \lab{1free_N} \ee
Again, by induction, it is easy to get
\be
\Phi_n(z)= \frac{1}{n+1} \sum_{k=0}^n (k+1)z^k, \quad n=N-1,N-2, \dots, 0. \lab{1free_Phi} \ee
Corresponding Verblunsky parameters are 
\be
a_n= -\frac{1}{n+2}, \; n=0,1,2,\dots, N-1, \; a_N =-1 \lab{a_single} \ee
This case corresponds to so-called "single momentum" POPUC \cite{Simon}.

They are orthogonal on vertices of regular $N+2$-gon apart from the vertex $z=1$. Corresponding concentrated masses can be easily calculated from \re{ws_st}
\be
w_s = (\zeta_s-1)(\zeta_s^{-1}-1)= 4 \sin^2 \theta_s, \quad s=1,2, \dots, N+1 \lab{w_s2} \ee 
where
\be
\theta_s = \frac{\pi s}{N+2} \lab{theta_ss} \ee
Note that the weight distribution \re{w_s2} is obtained from the distribution  $w_s=const$ by Christoffel transformation which deletes the spectral point $z=1$ (see \cite{CMG}).

Before starting investigation of more complicated cases connected with Kronecker polynomials, let us mention two simple statements which will be useful in what follows.

\vspace{3mm}

\begin{pr}
Assume that $\Phi_{N+1}(z)$ is an arbitrary $N+1$-order polynomial having all simple roots on the unit circle $|z|=1$. Let $\Phi_n(z), \: n=0,1,\dots, N$ be the corresponding unique sequence of monic Sturm OPUC fixed by the condition $\Phi_N(z) = (N+1)^{-1}\Phi_{N+1}'(z)$. Let $a_0, a_1, \dots, a_{N}$ be the corresponding Verblunsky coefficients (note that necessary $a_N =- {\bar \Phi}_{N+1}(0)$). 

Consider the polynomial
$\t \Phi_{N+1}(z)=(-1)^{N+1} \Phi_{N+1}(-z)$. Let  ${\t \Phi}_n(z), \: n=0,1,\dots, N$ be the corresponding unique sequence of monic Sturm OPUC fixed by the condition ${\t \Phi}_N(z) = (N+1)^{-1}{\t \Phi}_{N+1}'(z)$ with the Verblunsky parameters $\t a_n$.

Then 
\be 
\t \Phi_n(z) = (-1)^n \Phi_n(-z)
\lab{tPhi-} \ee
and
\be
\t a_n = (-1)^{n+1} a_n
\lab{ta-} \ee
\end{pr}

\vspace{3mm}

\begin{pr}
Under the same conditions let $k$ be an arbitrary positive integer and $\t N = k(N+1)-1$. Assume that  $\t \Phi_{{\t N}+1}(z)= \Phi_{N+1}(z^k)$. 

Then corresponding Sturm sequence of OPUC is 
\be 
\t \Phi_{nk+j}(z) = z^j \Phi_n(z^k), \; j=0,1,\dots, k-1, \; n=0,1,\dots, [N/k]
\lab{tPhik} \ee
and the Verblunsky parameters are
\be
\t a_{n} =   a_{j-1},\quad  \mbox{if} \quad n=kj-1, \; j=1,2,\dots 
\lab{tak} \ee
and $\t a_n=0$ otherwise.

\end{pr}  

Proof of these statements are elementary and can be done by induction and by uniqueness of OPUC defined by recurrence relation \re{rec_Phi}. Note that the OPUC with Verblusnky coefficients satisfying condition \re{tak} are called the {\it sieved} OPUC \cite{Ismail}, \cite{MS}.

The above propositions play the crucial role in the next section in constructing explicit examples of OPUC corresponding to a large family of cyclotomic polynomials.

\section{Sturmian cyclotomic POPUC}
\setcounter{equation}{0}
In this section we consider a special case of the Sturmian Kronecker polynomials when $K(z) = C_M(z)$, where $C_M(z)$ is a cyclotomic polynomial corresponding to the primitive $M$-th roots of unity. Recall that the degree of the polynomial $C_M(z)$ is 
\be
N+1=\varphi(M)
\lab{N+1_C} \ee

\vspace{3mm}

The orthogonality weights are (to within a normalization factor)
\be
w_s = \frac{1}{C_N'(\zeta_s) C_N'(\zeta_s^{-1})} = \frac{1}{|C_N'(\zeta_s)|^2}, \quad s=1,2 \dots d,
\lab{w_cyclo} \ee
where $\zeta_s$ are primitive roots of unity.

We consider first two simplest families of cyclotomic OPUC: when $N=p$ and $N=2p$, where $p$ is an odd prime number.

If $N=p$ then all $p$-th roots of unity (apart from 1) are primitive and hence 
\be
C_p(z) = \frac{z^p-1}{z-1} = z^{p-1} + z^{p-2} + \dots + z +1 \lab{C_p} \ee   
From the results of the previous section we have that in this case the Verblunsky parameters are
\be
a_n = -\frac{1}{n+2}, \: n=0,1,\dots, p-3, \: a_{p-2} =-1 
\lab{a_n_p} \ee

The case $N=2p$ is also simple. Indeed, in this case it is well known that \cite{Waerden}
\be
C_{2p}(z) =C_p(-z) \lab{C_2p} \ee
and that $N+1= \varphi(2p) = \varphi(p) = p-1$.

Using results of the previous section, we have the expression of the Verblunsky coefficients
\be
a_n = \frac{(-1)^n }{n+2}, \: n=0,1,\dots, p-3, \: a_{p-2} =-1 
\lab{a_2p} \ee
Corresponding OPUC are given by \re{tPhi-}.

Starting from these two simple examples, we can construct further explicit families by using the last two propositions of the previous section.  

Consider, e.g. the case $M=p^m$ where $p$ is a prime number. It is well known \cite{HW} that 
\be
C_M(z) = C_p\left(z^{p^{m-1}}\right)
\lab{C_M_m} \ee

Then from Proposition {\bf 4} it follows that corresponding POPUC $\Phi_n(z)$ are sieved polynomials given by \re{tPhik}.

Consider for example the case  $M=2^m, \: m=1,2,\dots$. The cyclotomic polynomial is
\be
C_{M}=z^{M/2} +1 =z ^{2^{m-1}} +1
\lab{C_M_2} \ee
Then all Verblunsky parameters are zero apart from the last one:
\be
a_0=a_1=\dots = a_{M-2}=0, \: a_{M-1}=-1.
\lab{a_2} \ee
This example is similar to the case of the "free" POPUC \re{a_s}. The only difference is the value of the last coefficient $a_N$. 

In the case $M=3^m$ the cyclotomic polynomial is 
\be
C_M = z^{2M/3} + z^{M/3} +1
\lab{C_M_3} \ee 
There are only two nonzero Verblunsky parameters: $a_{M/3-1}=-1/2$ and $a_{M-1}=-1$. All other parameters $a_n$ are zero.

Other cases can  be considered similarly.

We thus have the 
\begin{pr}
If $M=2^k p^m$ where $p$ is an odd prime number, then the Sturmiam cyclotomic POPUC $\Phi_n(z)$ are sieved polynomials having the only nonzero Verblunsky coefficients  $\frac{(-1)^n }{n+2}, \: n=0,1,\dots, p-3$ and $a_N=-1$ interlaced by zeros.  
\end{pr} 
For example for $M=2^2\cdot 5=20$ one has $C_{20}(z) = z^8 -z^6 +z^4-z^2+1$ and the sequence of Verblunsky parameters is $(0,1/2,0,-1/3,0,1/4,0,-1)$.

Consider a more complicated case of cyclotomic polynomial $\Phi_{N+1}(z)= C_{K}(z)$,  where $K=pq$ with two prime numbers $2<p<q$.

The degree of cyclotomic polynomial $C_{pq}(z)$ is $N+1=(p-1)(q-1)$. Hence, starting with $C_{pq}(z)$, we obtain the sequence $\Phi_0(z), \Phi_1(z), \dots, \Phi_N(z)$ of POPUC. These polynomials are orthogonal on the set $\zeta_k$ of primitive roots of the number $pq$
\be
\zeta_k = \exp\left( \frac{2 \pi i y_k}{pq}\right) , \quad k=1,2,3, \dots, (p-1)(q-1),
\lab{zeta_pq} \ee
where $y_k=1,2,\dots$ is the set of positive integers which are coprime with $pq$ and ordered by increasing. For example, when $p=3,q=5$ one has
\be
y_1=1,y_2=2,y_3=4,y_4=7,y_5=8,y_6=11,y_7=13,y_8=14.
\lab{y_15} \ee
The corresponding concentrated masses $w_k$ can be easily found from \re{ws_st}. Indeed, the cyclotomic polynomial $C_{pq}(z)$ for the two distinct prime numbers $p,q$ has the explicit expression \cite{HW}
\be
C_{pq}(z) = \frac{(z^{pq}-1)(z-1)}{(z^p-1)(z^q-1)}
\lab{C_pq_expl} \ee
Taking derivative of \re{C_pq_expl} and putting $z=\zeta_k$ we have from \re{ws_st} that the non-normalized weights are
 \be
w_k = \frac{\sin^2\left(\frac{\pi y_k}{p} \right)  \sin^2\left(\frac{\pi y_k}{q} \right)}{\sin^2\left( \frac{\pi y_k}{pq}\right)}, \quad k=1,2,3,\dots, (p-1)(q-1)
\lab{w_k_pq} \ee
 
The explicit expression for the Verblunsky parameters $a_n$ and for corresponding polynomials $\Phi_n(z)$ is a more complicated problem. We display here first nontrivial cases.

For $K=15=3 \cdot 5$ the cyclotomic polynomial is
\be
C_{15}(z) = {z}^{8}-{z}^{7}+{z}^{5}-{z}^{4}+{z}^{3}-z+1
\lab{cycl_15} \ee
while sequence of the Verblunsky parameters $a_i, \: i=0,1,\dots, N-1$ is
\be
\left\{ {2\over{3}} ,-{1\over{5}}, -{9\over{16}}, {1\over{5}}, -{2\over{7}}, {1\over{9}}, {1\over{8}} , -1\right\} 
\lab{a_15} \ee
For $K=21=3 \cdot 7$ the cyclotomic polynomial $C_{21}(z)$  is
\be
C_{21}(z) ={z}^{12}-{z}^{11}+{z}^{9}-{z}^{8}+{z}^{6}-{z}^{4}+{z}^{3}-z+1
\lab{cycl_21} \ee
and the sequence of Verblunsky parameters is 
\be
\left\{ {2\over{3}} ,-{1\over{5}}, -{1\over{4}}, {1\over{3}}, -{1\over{2}}, -{1\over{4}}, {1\over{10}}, {1\over{9}}, -{2\over{11}}, {1\over{13}} , {1\over{12}} , -1\right\} 
\lab{a_21} \ee

For generic case of $K=pq$ we have the following 

{\bf Conjecture}. Let $2<p<q$ be two prime numbers and let $n=pm+r$, where the remainder can take the values $r=0,1,2, \dots, p-1$.
Then for the "head" of the sequence $a_1, a_2, \dots, a_N$ with $0<n<q-p$ one has explicit expression of the Verblunsky parameters for corresponding OPUC 
\be
a_{pm+r} = \left\{ { \frac{p-1}{(m+1)p}, \quad \mbox{if} \quad r=0 \atop   -\frac{1}{(m+2)p-r} \quad \mbox{otherwise}}  \right.
\lab{a_n<q-p} \ee
For the "tail" of this sequence we have the expression (which is valid for $n \le q-2$)
\be
a_{N-n-1} = \left\{ { -\frac{p-1}{N-n-1+p}, \quad \mbox{if} \quad p|(n+1) \atop   \frac{1}{N+1-n+2\cdot mod(n,p)} \quad \mbox{otherwise}}  \right.,
\lab{a_n<p-2} \ee
where $N=(p-1)(q-1)-1$. Observe that the structure of both "head" and "tail" depends entirely on the smaller prime factor $p$. The bigger prime factor $q$ only regulates the "length" of the sequences where rules \re{a_n<q-p} and \re{a_n<p-2} are valid .

This conjecture was verified (with the help of computer) for all pairs of primes $p<q$ until $q=29$. 

As an illustration of the above conjecture we display hear the "head" and the "tail" of the Verblunsky parameters for the case $p=5, q=17$. In this case $N=63$. The first parameters $a_0, a_1, \dots, a_{12}$ are:
\be
\left\{ \frac{4}{5},-\frac{1}{9},-\frac{1}{8},-\frac{1}{7},-\frac{1}{6}, \frac{2}{5},-\frac{1}{14},-\frac{1}{13},-\frac{1}{12},-\frac{1}{11},{\frac {4}{15}},-
\frac{1}{19},-{\frac {65}{144}}, \dots \right\}
\ee
The last parameters $a_{46}, a_{47}, \dots, a_{N}$ are:
\be
\left\{\dots, -{\frac {481}{1920}},\frac{1}{49},-{\frac {4}{53}},{\frac {1}{57}},{\frac {1}
{56}},{\frac {1}{55}},{\frac {1}{54}},-{\frac {2}{29}},{\frac {1}{62}}
,{\frac {1}{61}},{\frac {1}{60}},{\frac {1}{59}},-{\frac {4}{63}},{
\frac {1}{67}},{\frac {1}{66}},{\frac {1}{65}},{\frac {1}{64}},-1 \right\}
\ee

In the middle part the sequence $[a_0,a_1, \dots, a_N]$ for the cyclotomic polynomial $C_{pq}(z)$ becomes more complicated, for example for the case $p=5,q=7$ the first 14 members of this sequence are
\be
\left\{\frac{4}{5},-\frac{1}{9},-{\frac {13}{32}},-{\frac {5}{19}},-{\frac {4}{15}},{\frac {
25}{77}},-{\frac {139}{884}},{\frac {1049}{4619}},-{\frac {302}{1635}}
,-{\frac {204}{559}},{\frac {359}{2982}},-{\frac {1292}{15677}},{
\frac {20566}{163715}},-{\frac {6099}{29521}},\dots \right\}
\lab{5*7} \ee

Finally, we present explicit formulas for two last Verblunsky coefficients $a_{N}$ and $a_{N-1}$ for any Sturmian cyclotomic polynomial.

\begin{pr}
Let $M>1$ be an arbitrary positive integer.  Put $N=\varphi(M)-1$ and $\Phi_{N+1}(z) = C_M(z)$. Let $\Phi_{N}(z), \Phi_{N-1}(z), \dots, \Phi_0(z)=1$ be corresponding Sturmian chain of POPUC. Then $a_N = -1$ and $a_{N-1}=\frac{\mu(M)}{N+1}$, where $\mu(n)$ is the M\"obius function.
\end{pr}
{\it Remark}. Recall that the M\"obius function $\mu(n)$ is defined on positive integers $n$ as $(-1)^j$ if the decomposition of $n$ into $j$ prime factors is square free and $\mu(n)=0$ if this decomposition is not square free. 

The proof of this proposition follows from the definition of Verblunsky coefficients. Namely
\be 
a_N = -\Phi_{N+1}(0) =-1
\lab{a_N_any} \ee  
because the free term of any cyclotomic polynomial $C_M(z)$ with $M>1$ is equal to 1 \cite{HW}. Similarly, $a_{N-1}=-\frac{\Phi_{N+1}'(0)}{N+1}$. But  $\Phi_{N+1}'(0)$ is equal to the coefficient in front of the linear term $z$ of the cyclotomic polynomial $C_M(z)$. From the Vieta theorem and from the palindromic property $z^{N+1} \Phi_{N+1}(1/z) = \Phi_{N+1}(z)$ it follows that this coefficient is equal (up to a sign) to sum of all primitive roots $\zeta_k$ of unity. And this sum is known to be $\mu(M)$ \cite{HW}
\be
-\Phi_{N+1}'(0) =  \sum_{k=1}^M \zeta_k = \mu(M). \lab{sum_M} \ee

\section{Beyond cyclotomic polynomials}
\setcounter{equation}{0}
The anti-cyclotomic (or inverse cyclotomic \cite{HM}) polynomial $A_N(z)$  were introduced in \cite{Moree} as the minimal polynomial for all {\it non-primitive} $N$-th roots of unity. In other words all roots of the anti-cyclotomic polynomial $A_N(z)$ are complementary to the roots of $C_N(z)$ among all $N$-th  roots of unity.  This leads to the explicit formula \cite{HM}
\be
A_N(z) = \frac{z^N-1}{C_{N}(z)}. \lab{A_N_exp} \ee  
Clearly $d=\deg(A_N(z)) = N-\vphi(N)$.

Of course, anti-cyclotomic polynomials belong to special types of Kronecker polynomials. Hence, one can consider corresponding Sturmian OPUC sequences.

Note that the anti-cyclotomic polynomials are {\it anti-palindromic}, i.e.  they have the property
\be
z^d A_N(1/z) = -A_N(z).
\lab{anti-pal} \ee
This property is obvious because $z=1$ is always a root of any anti-palyndromic polynomial.

In this section we consider only the simplest nontrivial case when $N=2p$ where $p$ is a prime number $p \ge 3$. Then
\be
A_{2p}(z)=  z^{p+1} + z^{p} -z-1
\lab{A_{2p}_exp} \ee
By induction, it is easy to check that the Verblunsky parameters are 
\be 
a_0 = \frac{1-p}{p}, \quad a_n = \frac{(-1)^{n}}{2p-n}, \; n=1, 2, \dots, p-1, \quad a_{p}=1.
\lab{a_n_A} \ee
Corresponding POPUC have the explicit expression:
\be
\Phi_0(z)=1, \: \Phi_1(z) = z+ 1-p^{-1}, \quad \Phi_n(z) = {z}^{n}+{\frac { \left( 2\,p-n \right) {z}^{n-1}}{2\,p-n+1}}+{\frac {
 \left( -1 \right) ^{n}}{2\,p-n+1}}, \; n=2,3,\dots, p.
\lab{Phi_anti} \ee 
Expression \re{Phi_anti} can easily be checked by recurrence relation \re{rec_Phi}. Note that any polynomial $\Phi_n(z)$ for $n\ge 2$ contains only three terms: $z^n, z^{n-1}$ and $z^0$.

The orthogonality relation for these anti-cyclotomic polynomials can easily be found from general formula \re{ws_st}. Indeed, the only spectral points are vertices of the regular $p$-gon
\be
\zeta_k = \exp\left( \frac{2 \pi i k}{p}\right), \quad k=0,1,2,\dots, p-1  
\lab{zeta_p} \ee
with the concentrated masses
\be
w_k = \frac{1}{2p^2 \cos^2\left( \frac{\pi  k}{p}\right) }
\lab{w_k_anti} \ee
and the additional point $z=-1$ with the concentrated mass $w_{-1}=1/2$. So we have the orthogonality relation
\be
\sum_{k=0}^{p-1} \frac{\Phi_n(\zeta_k) \Phi_m(\zeta_k^{-1})}{2p^2 \cos^2\left( \frac{\pi  k}{p}\right)}  + \frac{\Phi_n(-1) \Phi_m(-1)}{2} = h_n \delta_{nm}
\lab{ort_anti} \ee
Note this measure is appropriately normalized, i.e.
\be
\sum_{k=0}^{p-1} w_k + w_{-1} =1.
\lab{norm_anti} \ee
Formula \re{norm_anti} follows from the elementary trigonometric identity
\be
\sum_{k=0}^{M-1} \sec^2\left( \frac{\pi  k}{M}\right) =M^2 \lab{elem_sec} \ee
which is valid for any positive odd number $M$.

Finally, consider only one simple (but leading to nontrivial results) example of Kronecker polynomials beyond the cyclotomic or anti-cyclotomic polynomials. 

Namely, we {\it adjoin} one additional root $z=-1$ to the  polynomial $\frac{z^M-1}{z-1}$ for an {\it odd} (not necessarily prime) $M$:
\be
K_M(z) = (z+1)(z^{M-1} + z^{M-2} + \dots + z+1) 
\lab{K_M} \ee 
The degree of this polynomial is $M$. Hence $N=M-1$ and we have for the polynomial $\Phi_{N}(z)$
\be
\Phi_N(z)=\Phi_{M-1}(z) = M^{-1} \left(M z^{M-1} + 2(M-1)z^{M-2} + 2(M-2)z^{M-3} + \dots + 6 z^2 + 4 z + 2 \right) 
\lab{Phi_M-1} \ee
Then it is possible to show that the Verblunsky parameters for $n=0,1,\dots, M-2$ are given by
\be
a_n = \left\{ {-\frac{2(n-M+1)}{(n+1)(n-2M+1)}, \quad   \mbox{even  n} \atop \frac{2M}{n^2  - 2(M-1)n -2M} , \quad   \mbox{odd  n}}  \right.
 \lab{a_nca} \ee
and
\be
a_{M-1}=-1
\lab{a_M-1} \ee
The polynomials $\Phi_n(z)$ have simple explicit expression 
\be
\Phi_n(z) =z^n +  \sum_{k=0}^{n-1} \left( -a_{n-1} + \beta_n k  \right) z^k, 
\lab{expl_Phi_adj} \ee
where
\be
\beta_n = \left\{   {{\frac {2(2\,M-n-1)}{2\,Mn-{n}^{2}+1}}, \quad n \; \mbox{even} \atop \frac{2}{n}, \quad n \; \mbox{odd} }   \right.
\lab{beta_n} \ee
Again, formulas \re{a_nca} and \re{expl_Phi_adj} can be verified directly by using recurrence relation \re{rec_Phi}.

It is interesting to note that the expansion coefficients $A_{nk}$ of the polynomial $\Phi_n(z)= z^n + \sum_{k=0}^{n-1} A_{nk}z^k$ depend on $k$ {\it linearly}.

Using formula \re{w_s_OPUC} one can find orthogonality measure for the obtained polynomials. Indeed, the only points of concentrated masses on the unit circle are zeros of the polynomial  $K_M$, i.e. there are $M-1$ points 
\be
q_k = \exp\left(\frac{2 \pi i k}{M}\right), \quad k=1,2,\dots, M-1 \lab{q_k} \ee
with corresponding weights 
\be
w_k = M^{-2} \tan^2\left( \frac{k \pi}{M}\right), \lab{w_k} \ee
and the additional point at $z=-1$ with the weight
\be
w_{-1}=1 \lab{w-1} \ee
Hence the orthogonality relation takes the form
\be
\sum_{k=1}^{M-1} M^{-2} \tan^2\left( \frac{k \pi}{M}\right) \Phi(q_k,n)\Phi(1/q_k,m) + \Phi(-1,n)\Phi(-1,m) = h_n \delta_{nm} \lab{ort_M} \ee
The normalization constant
\be
h_0 = 2-M^{-1}. \lab{h_0_M} \ee
Formula \re{h_0_M} follows from the elementary trigonometric identity
\be
\sum_{k=1}^{M-1} \tan^2\left(\frac{k \pi}{M} \right) = M(M-1) \lab{sum_tan} 
\ee
which is valid for odd $M$. 

\vspace{3mm}

{\it Remark}. Condition that $M$ is {\it odd} in \re{K_M} is necessary. Indeed, if $M$ is even that the Kronecker polynomial $K_M(z)$ has double root at $z=-1$. This contradicts to the condition that all roots of the polynomial $K_M'(z)$ should lie {\it inside} the unit circle.

\vspace{5mm}

\section{Conclusions}
\setcounter{equation}{0}
We have demonstrated that numerous examples of the Kronecker polynomials give rise to new explicit finite families of POPUC with the measure located at roots of unity. 

There are many remaining problems for further studying. One of them is finding explicit expression of the polynomials $\Phi_n(z)$ when the Kronecker polynomial coincide with either binary $C_{pq}(z)$ or with ternary $C_{pqr}(z)$ cyclotomic polynomials, where $p,q,r$ are distinct odd prime numbers.  Another problem is studying possible infinite families in the limiting cases (say, when $ p \to \infty$ in \re{w_k_pq}). It is expected to lead to some nontrivial examples of POPUC orthogonal with respect to a continuous measure.

\vspace{6mm}

{\large \bf Acknowledgments.} 

\vspace{5mm}

 The author is gratefully holding Simons CRM Professorship. The work is funded by the National Foundation of China (Grant No.11771015).

\bb{99}

\bi{CMG} K.Castillo, F.Marcellan and L. Garcia, {\it Linear spectral transformations, Hessenberg matrices, and orthogonal
polynomials}, Rend. Circ. Matem. Palermo, 
Serie II, Suppl. {\bf 82} (2010), 1--24.

\bi{CMR} K.~Castillo, F.~Marcell\'an and M.~N.~Rebocho    {\it Zeros of para-orthogonal polynomials and linear spectral transformations on the unit circle}, Numer. Algor. {\bf 71}, 699--714(2016).

\bi{Chi} T. Chihara, {\it An Introduction to Orthogonal Polynomials},
Gordon and Breach, NY, 1978.

\bi{Dam} P.A.Damianou, {\it Monic Polynomials in $Z[x]$ with Roots in the Unit Disc}, Am.Math.Month., {\bf 108} (2001), 253--257.

\bi{Ger} Ya.L. Geronimus,\quad {\it Polynomials Orthogonal on a
Circle and their Applications}, \\ Am.Math.Transl.,Ser.1 {\bf
3}(1962), 1-78.


\bi{HW} G.Hardy and E.M.Wright, An introduction to the theory of numbers (Sixth ed.), 2008, Oxford University Press.

\bi{HM}  A. Herrera-Poyatos, Pieter Moree, {\it Coefficients and higher order derivatives of cyclotomic polynomials: Old and new}, Expos. Math. (to appear), ArXiv: 1805.05207.

\bi{Ismail} M.E.H.Ismail and Xin Li, {\it On sieved orthogonal polynomials IX: Orthogonality on the unit circle}, Pacific J.Math. {\bf 153}, (1992), 289--297.


\bi{Kron} L.~Kronecker, {\it Zwei S\"atze 
\"uber Gleichungen mit ganzzahligen Coefficienten}, Crelle, Oeuvres
I (1857) 105--108.

\bi{MS} F.Marcellan, A. Sri Ranga, {\it Sieved para-orthogonal polynomials on the unit circle}, Appl. Math. and Comp.,  {\bf 244} (2014) 335--343.

\bi{Mar1} A.~Mart\'inez-Finkelshtein, B.~Simanek, and B.~Simon. {\it Poncelet's theorem, paraorthogonal polynomials
and the numerical range of compressed multiplication operators}. Adv. Math., {\bf 349}, 992--1035, (2019).

\bi{Moree} P.Moree, {\it Inverse cyclotomic polynomials}, J.Numb.Th. {\bf 129} (2009) 667--680.

\bi{Simon} B.Simon, {\it Orthogonal Polynomials On The Unit
Circle}, AMS, 2005.

\bi{Waerden} B.Van der Waerden, Algebra, vol. I, Springer, 1991.

\bi{ZheS} A.Zhedanov, {\it Classical Sturmian sequences}, arXiv:1904.03789.


\eb

\end{document}